\documentstyle{amsppt}
\magnification=1200
\hsize=6.5truein
\vsize=8.9truein
\topmatter
\title Exponents of an irreducible plane curve singularity
\endtitle
\author   Morihiko Saito
\endauthor
\affil RIMS Kyoto University, Kyoto 606-8502 Japan \endaffil
\keywords exponents, Puiseux pairs, Enriques diagram, mixed Hodge structure
\endkeywords
\subjclass 32S35, 32S40\endsubjclass
\endtopmatter
\tolerance=1000
\baselineskip=12pt 
\def\scirc{\raise.2ex\hbox{${\scriptstyle\circ}$}}
\def\ssbull{\raise.2ex\hbox{${\scriptscriptstyle\bullet}$}}
\def\Gr{\text{{\rm Gr}}}
\def\Ker{\hbox{{\rm Ker}}}

\def\cond[#1]{\par\noindent\rlap{(#1)}\hskip1.5cm
\hangindent=1.5cm\hangafter=1}

\document

\noindent
Let 
$ f : ({\Bbb C}^{n+1},0) \rightarrow  ({\Bbb C},0) $ be a germ of a 
holomorphic function with an isolated singularity.
Using Steenbrink's theory [15] of mixed Hodge structure on 
the cohomology of the Milnor fiber, we can define the {\it exponents} 
(or spectra, up to the shift by one, in the terminology of Varchenko 
[16]) to be
$ \mu  $ rational numbers
$ \{\alpha _{1}, \dots , \alpha _{\mu }\} $ such that
$ \exp(2\pi i(\sqrt {-1}\alpha _{i}) $ are the eigenvalues of the Milnor 
monodromy and their integral part is determined by the Hodge filtration 
of the mixed Hodge structure.
This notion was first introduced by Steenbrink [15].
 It is known that the exponents are constant under
$ \mu  $-constant deformation of
$ f $.
See [18].
In particular, they depend only on
$ f^{-1}(0) $.
They express the vanishing order (up to the shift by one) of the period 
integrals of holomorphic forms on vanishing cycles.
See [16], [17].

Let
$ (V,0) \subset  ({\Bbb C}^{2}, 0) $ be a germ of a reduced and 
irreducible plane curve defined by a holomorphic function
$ f $.
It is known that the equisingular class of
$ (V,0) $ is determined by its {\it Puiseux pairs}, and numerical 
invariants such as the Milnor number or the characteristic polynomial 
of the monodromy can be expressed in terms of the Puiseux pairs, 
cf. [2], [8].
In this note we give an explicit formula for the exponents of
$ f $ in terms of the Puiseux pairs, cf. Theorem (1.5).
The proof uses the {\it Enriques diagram} [3], [4] of
$ (V,0) $ which describes the canonical process of embedded resolution 
of a plane curve
$ V $ by iterating point center blowing-ups.
In the irreducible case, we can describe explicitly the Enriques diagram
as well as the multiplicity of the pull-back of
$ f $ along the irreducible components of the exceptional divisor of the
resolution, using the continued fraction expansion of the Puiseux pairs.
Then we can apply a formula of Steenbrink [15] to calculate the Hodge 
numbers of the vanishing cohomology.

As an application we can prove in this case a recent conjecture of 
Hertling [7] that the variance (i.e. the square of the standard deviation) 
of the exponents is bounded by the difference between the maximal 
and minimal exponents divided by 
$ 12 $.  
See (5.2).  
This was rather unexpected, because no philosophical reason for the 
conjecture is known.  
For the proof, we need a good estimate of the average of the exponents 
less than  
$1 $ in the quasihomogeneous case  (see (5.3)) which is sharp and leads 
to a simple expression like (5.2.3).
The proof shows that his estimate of the variance is rather sharp
in our case.  
We might expect that his conjecture should hold in a more general case.

I would like to thank D.T. L\^e, K. Saito, 
J. Steenbrink and T. Yano for useful discussions and for
their interest in this problem, and C. Hertling for verifying the argument 
and informing me of nontrivial misprints.
I also thank A. Dimca for drawing my attention to Hertling's 
conjecture.
This paper (except for Sect. 5) is a revised version of a paper typed 
in 1990 at RIMS, and its original manuscript (which is quoted in [13]) 
was written in 1982.

\bigskip\bigskip
\centerline{\bf 1. Exponents}

\bigskip
\noindent
{\bf 1.1.}
Let
$ f : ({\Bbb C}^{n+1},0) \rightarrow  ({\Bbb C},0) $ be a germ of a holomorphic
function with an isolated singularity.
Let
$ f : X \rightarrow  S $ be a good representative of
$ f  $ (sometimes called a Milnor fibration) defined by
$$
X = \{x \in  {\Bbb C}^{n+1} : |x| < \varepsilon , |f(x)| < \eta \}, \quad
S = \{t \in  {\Bbb C} : |t| < \eta \}
$$
for 
$ 0 < \eta  \ll \varepsilon  \ll 1 $.

As in [15], we denote by
$ H^{n}(X_{\infty },{\Bbb C}) $ the vanishing cohomology of
$ f,  $ which is (noncanonically) isomorphic to the cohomology of the
Milnor fiber
$ H^{n}(X_{t},{\Bbb C}) $ for
$ t \ne  0 $.
In [15] Steenbrink constructed a canonical mixed Hodge structure on
$ H^{n}(X_{\infty },{\Bbb C}) $ using Deligne's theory of mixed Hodge
structure [5].
Let
$ F $ be the Hodge filtration on
$ H^{n}(X_{\infty },{\Bbb C}) $.
Let
$ \mu  = \dim_{{\Bbb C}}H^{n}(X_{\infty },{\Bbb C}),  $ which is called
the Milnor number of
$ f $.
Let
$$
H^{n}(X_{\infty },{\Bbb C})_{\lambda } = \Ker(T_{s} - \lambda  : 
H^{n}(X_{\infty },{\Bbb C}) \rightarrow  H^{n}(X_{\infty },{\Bbb C})),
$$
where
$ T = T_{s}T_{u} $ is the Jordan decomposition of the monodromy.

\medskip
\noindent
{\bf 1.2.~Definition} [15].
The {\it exponents} of
$ f $ are
$ \mu  $-rational numbers
$ \{\alpha _{1}, \dots , \alpha _{\mu }\} $ such that
$$
0 < \alpha _{1} \le  \alpha _{2} \le  \dots \le  \alpha _{\mu } 
<n+1
\tag 1.2.1
$$
and are defined by the following condition:
$$
\aligned
\#\{j : \bold{e}(-\alpha _{j}) = \lambda , [\alpha _{j}] = n - p\}
&= \dim_{{\Bbb C}}{\Gr}_{F}^{p}
H^{n}(X_{\infty },{\Bbb C})_{\lambda } \,(\lambda  \ne  1)
\\
\#\{j : \alpha _{j} = n - p + 1\}   
&= \dim_{{\Bbb C}}{\Gr}_{F}^{p}H^{n}(X_{\infty },{\Bbb C})_{1}
\endaligned
\tag 1.2.2
$$
where
$ [\alpha ] = \max\{i \in  {\Bbb Z} : i \le  \alpha \} $ is the Gauss symbol, 
$ \bold{e}(\alpha ) = \exp(2\pi i\alpha ),  $ and
$ F $ is the Hodge filtration of the mixed Hodge structure on
$ H^{n}(X_{\infty },{\Bbb C}) $.
(See also [10].)
Let
$$
\chi _{f}(t) = \sum _{1\le i\le \mu } t^{\alpha _{i}}
$$
which is called the spectrum of
$ f $ by a recent terminology.

\medskip\noindent
{\it Remark.}
By the symmetry of the exponents (cf. [15]), we have
$$
\chi _{f}(t) = t^{n+1}\chi _{f}(t^{-1}).
\leqno(1.2.3)
$$

\noindent
{\bf 1.3.}~{\it Example} (quasihomogeneous case).
Let
$ f $ be a quasihomogeneous polynomial with weight
$ (w_{0}, \dots , w_{n}) $, i.e.,
$ f $ is a linear combination of monomials
$ {x}_{0}^{{m}_{0}}{\cdots}{x}_{n}^{{m}_{n}} $ such that
$ \sum  w_{i}m_{i} = 1 $.
Then
$$
\Omega _{f} := {\Omega }_{X,0}^{n+1} / df\wedge {\Omega }_{X,0}^{n}
\leqno(1.3.1)
$$
is a graded 
$ \mu  $-dimensional vector space whose grading is induced by the weight
of the coordinates and is indexed by
$ \bold{Q} $.
Let
$ P(\Omega _{f},t) $ be the Poincar\'e polynomial of
$ \Omega _{f} $.
By [14], we have
$$
P(\Omega _{f},t) = \chi _{f}(t).
\leqno(1.3.2)
$$
Furthermore
$ P(\Omega _{f},t) $ is written explicitly in terms of the weights, 
and we get
$$
\chi _{f}(t) = \prod\nolimits\limits_{i=0}^{n}{{t}^{{w}_{i}}-t \over 1-{t}^{{w}_{i}}}
\leqno(1.3.3)
$$
as is well-known.
This formula is generalized to the nondegenerate Newton polyhedron case
by [15], [11] using the Newton filtration on
$ \Omega _{f} $.

\medskip
\noindent
{\bf 1.4.~Definition.} Let
$ (k_{i}, n_{i}) $ be pairs of relatively prime positive integers such that
$ n_{i} > 1 $ for 
$ 1 \le  i \le  g $.
Let
$ w_{i} $ be integers defined inductively by
$ w_{1} = k_{1} $ and
$$
w_{i} = w_{i-1}n_{i-1}n_{i} + k_{i} \quad\text{for $ i > 1 $.}
$$
We define
$ \Phi _{g}(k_{1},n_{1}; \dots ;k_{g},n_{g})(t) $ by induction on
$ g $ as follows:
$$
\alignat2
\Phi _{1}(k,n)(t) = 
&\left({{{t}^{1/k}-t \over 1-{t}^{1/k}}}\right) 
\left({{{t}^{1/n}-t \over 1-{t}^{1/n}}}\right)
&&\quad \text{if }g = 1,
\tag 1.4.1
\\
\Phi _{g}(k_{1},n_{1}; \dots ;k_{g},n_{g})(t) = 
&\left({{{t}^{1/{w}_{g}}-t \over 1-{t}^{1/{w}_{g}}}}\right)
\left({{{t}^{1/{n}_{g}}-t \over 1-{t}^{1/{n}_{g}}}}\right) 
&&
\tag 1.4.2
\\
+
&\left({{1-t \over 1-{t}^{1/{n}_{g}}}}\right)
{\Phi }_{g-1}^{<1}(t^{1/n_{g}}) 
&&
\\
+
&\,\,{t^{1-1/n_{g}}}\left({{1-t \over 1-{t}^{1/{n}_{g}}}}\right)
{\Phi }_{g-1}^{>1}(t^{1/n_{g}})
&&\quad\text{if }g > 1,
\endalignat
$$
where
$ {\Phi }_{g-1}^{<1}(t) $ and
$ {\Phi }_{g-1}^{>1}(t) $ are defined respectively by
$ \sum _{\alpha <1} c_{\alpha }t^{\alpha } $ and
$ \sum _{\alpha >1} c_{\alpha }t^{\alpha } $ with
$ \sum _{\alpha } c_{\alpha }t^{\alpha } = \Phi _{g-1}(k_{1},
n_{1}; \dots ;k_{g-1},n_{g-1})(t) $

\proclaim{{\bf 1.5.~Theorem}}
Let 
$ (V,0) \subset  ({\Bbb C}^{2},0) $ be a reduced and irreducible curve defined
by a function
$ f,  $ and
$ (k_{1}, n_{1}), \dots , (k_{g}, n_{g}) $ the Puiseux pairs of
$ V $, cf. (2.1).
Then
$$
\chi _{f}(t) = \Phi _{g}(k_{1},n_{1}; \dots ;k_{g},n_{g})(t).
\leqno(1.5.1)
$$
\endproclaim

\noindent
{\it Remark.}
The assertion means that the exponents which are smaller than 
$ 1 $ are given (with multiplicity) by
$$
\left\{{\left({{1 \over {n}_{\nu +1}\cdots
{n}_{g}}\left({{i \over {n}_{\nu }}+{j \over {w}_{\nu }}}\right)+
{r \over {n}_{\nu +1}\cdots {n}_{g}}}\right)}\right\}
\leqno(1.5.2)
$$
for
$ 0 < i < n_{\nu } $,
$ 0 < j < w_{\nu } $,
$ 0 \le  r < n_{\nu +1} \cdots n_{g} $,
$ 1 \le \nu \le g $ such that
$ i/n_{\nu } + j/w_{\nu } < 1 $.

\bigskip\bigskip
\centerline{\bf 2. Puiseux pairs and Enriques diagram}

\bigskip
\noindent
{\bf 2.1.~Definition.} Let
$ (V,0) $ be a germ of a reduced and irreducible plane curve in
$ ({\Bbb C}^{2},0) $.
Let
$ (x,y) $ be the coordinate system of
$ ({\Bbb C}^{2},0) $.
We have the Puiseux expansion associated with
$ V $:
$$
\aligned
y = \,&\sum _{1\le i\le [k_{1}/n_{1}]} c_{0,i}x^{i}
\\
+\,&\sum _{0\le i\le [k_{2}/n_{2}]} c_{1,i}x^{(k_{1}+i)/n_{1}}
\\
+\,&\sum _{0\le i\le [k_{3}/n_{3}]} c_{2,i}x^{k_{1}/n_{1} + 
(k_{2}+i)/n_{1}n_{2 }  } 
\\
&\qquad\qquad\vdots
\\
+\,&\quad\,\,\,\sum _{0\le i} c_{g,i}x^{k_{1}/n_{1} + k_{2}/n_{1}n_{2} + {\cdots} + 
(k_{g}+i)/n_{1}{\cdots}n_{g}}
\endaligned
\leqno(2.1.1)
$$
where
$ c_{j,i} \in  {\Bbb C} $ such that
$ c_{j,0} \ne  0 \,(j \ne  0),  $ and
$ k_{j}, n_{j} \in  {\Bbb Z}_{+} $ such that
$ (k_{j}, n_{j}) = 1 $ and
$ n_{j} > 1 $.
In this paper, we call
$ (k_{1}, n_{1}), \dots , (k_{g}, n_{g}) $ the {\it Puiseux pairs} of
$ (V,0) $ with respect to the coordinates
$ x, y $.
It is known that the Puiseux pairs are independent of the coordinates as long
as the condition
$ k_{1} > n_{1} $ is satisfied (cf. for example [1], [19]).
This follows also from (2.6) below.
We will assume always this condition, unless the coordinates are specified
explicitly.
Note that the condition is always satisfied by exchanging the coordinates
$ x, y $ if necessary.

We define the {\it modified Puiseux pairs} of
$ (V,0) $ with respect to the coordinates
$ x, y $ by deleting
$ \sum _{1\le i\le [k_{1}/n_{1}]} c_{0,i}x^{i} $ in (2.1.1) and allowing
$ n_{1} = 1 $ if
$ c_{0,i} \ne  0 $ for some
$ i $ in the old expression.
If
$ n_{1} > 1 $ (i.e., 
$ c_{0,i} = 0 $ in the old expression), the modified 
Puiseux pairs are the same as the Puiseux pairs with respect to the 
coordinates.
(This notion depends on the coordinates, and will be used for the 
inversion of Puiseux pairs, cf. (2.7).)

\medskip\noindent
{\it Remark.} If we are interested only in the topological type, we may 
assume that
$ c_{j,i} = 0 $ for
$ i > 0 $ by deforming
$ V $ with the topological type unchanged.
However, we cannot do this if we consider a reducible curve.
For example, if it has two irreducible components defined by
$ y = x^{7/2} $ and
$ y = x^{a} + x^{7/2} $ with
$ a = 2 $ or
$ 3 $,
then the
$ EN $-diagram [6] consists of three splice components, and two of them are
$ \Sigma (1,2,7) $ (where
$ 1 $ corresponds to the proper transform of each irreducible 
component), but the middle splice component is
$ \Sigma (1,1,a) $.

\medskip\noindent
{\bf 2.2.}
Let
$ (V,0) $ be a germ of a reduced plane curve in
$ ({\Bbb C}^{2},0) $.
We have a canonical embedded resolution
$ \phi  : (X', D) \rightarrow  ({\Bbb C}^{2},0) $ by iterating point center 
blow-ups along the points at which the total transform of
$ V $ does not have normal crossings, where
$ D = \phi ^{-1}(0) $.
Let
$ V' $ be the proper transform of
$ V $.
By definition,
$ \phi ^{-1}(V) = D \cup  V' $ is a divisor with normal crossings, and the 
irreducible components
$ D_{\alpha } (\alpha  \in  \Lambda )  $ of
$ D $ are
$ \bold{P}^{1} $.
We say that
$ D_{\alpha } $ is proximate to
$ D_{\beta } $ if
$ D_{\beta } $ is first obtained as (the proper transform of) the exceptional
divisor
$ E_{\beta } $ of a blow-up and then
$ D_{\alpha } $ is obtained as (the proper transform of) the exceptional
divisor of a blow-up whose center is contained in the proper transform of
$ E_{\beta } $, cf. [3], [4], [19], [20].

\medskip
\noindent
{\bf 2.3.~Definition} (cf. [3], [4]).
With the above notation, the {\it Enriques diagram} of
$ V $ is an oriented graph
$ E $ consisting of white and black vertices and arrows:

\noindent
(i) The white vertices correspond to the exceptional divisors
$ D_{\alpha } \,(\alpha  \in  \Lambda ),  $ and the black vertices to
the irreducible components of
$ V' $.

\noindent
(ii) There is an arrow from a white vertex corresponding to
$ D_{\alpha } $ to a white vertex corresponding to
$ D_{\beta } $ if
$ D_{\alpha } $ is proximate to
$ D_{\beta } $.
There is an arrow form a black vertex to a white vertex if the corresponding
divisors intersect.
There are no other arrows.

\medskip
\noindent
{\bf 2.4.}~{\it Continued fraction expansion.}
Let
$ k, n $ be relatively prime positive integers such that
$ n > 1 $.
Let
$$
{k \over n} = a_{1} {+} {1 \over {a}_{2}} 
\raise-1.5ex\hbox{$+ \cdots  +$} 
{1 \over {a}_{h}}
\leqno(2.4.1)
$$
be the continued fraction expansion, i.e., there are positive integers
$ r_{0}, \dots , r_{h} $ such that
$ r_{0} = k, r_{1} = n $ and
$$
r_{i-1} = a_{i}r_{i} + r_{i+1},\quad 
a_{i} = \left[{r_{i-1}\over r_{i}}\right]\quad  \text{for }1 
\le  i \le  h,
\leqno(2.4.2)
$$
with
$ r_{h+1} = 0,  $.
Here
$ r_{h} = 1 $ because
$ (k, n) = 1 $.
We have
$ a_{1} = 0 $ if
$ k < n $.

We define nonnegative integers
$ P_{1}, \dots , P_{h} $ and
$ Q_{1}, \dots , Q_{h} $ by
$$
\bold{v}_{i} = a_{i}\bold{v}_{i-1} + \bold{v}_{i-2}\,\,\, \text{for }
1 \le  i \le  h
\leqno(2.4.3)
$$
with
$ \bold{v}_{i} = (P_{i}, Q_{i}),  $ where
$ \bold{v}_{-1} = (0,1), \bold{v}_{0} = (1,0) $.
Then
$$
P_{i}Q_{i-1} - P_{i-1}Q_{i} = (-1)^{i},
\leqno(2.4.4)
$$
and
$ (\bold{v}_{i-1},\bold{v}_{i}) $ is a basis of
$ {\Bbb Z}^{2}  $.
So there exist integers
$ x_{i}, y_{i} \,(0 \le  i \le  h) $ such that
$$
(k,n) = x_{i}\bold{v}_{i} + y_{i}\bold{v}_{i-1}.
$$
Note that
$ x_{0} = k, y_{0} = n $ and
$ x_{i-1} = a_{i}x_{i} + y_{i}, y_{i-1} = x_{i} \,(1 \le  i \le  h) $ by 
(2.4.3).
Comparing this with (2.4.2), we get
$ x_{i} = r_{i} $ and
$ (x_{h}, y_{h}) = (1,0) $, i.e.,
$$
P_{h} = k, Q_{h} = n.
\leqno(2.4.5)
$$
This implies
$$
{{P}_{i} \over {Q}_{i}} = a_{1} {+} {1 \over {a}_{2}} 
\raise-1.5ex\hbox{$+ \cdots  +$} {1 \over {a}_{i}}.
\leqno(2.4.6)
$$
by replacing
$ k/n $ with the right hand side of (2.4.6).

\medskip
\noindent
{\bf 2.5.}
With the above notation, we define the oriented graph
$ E(k,n) $ as follows:

\noindent
(i) The vertices of
$ E(k,n) $ consist of white vertices
$ \{D_{i,j} : 1 \le  i \le  h, 1 \le  j \le  a_{j}\} $ and a black vertex
$ V' $.

\noindent
(ii) There is an arrow from
$ D_{i',j'} $ to
$ D_{i,j} $ if one of the following conditions is satisfied:

(a)
$ i' = i $ and
$ j' = j + 1 $,

(b)
$ i' = i + 1 $ and
$ j = a_{i} $,

(c)
$ i' = i + 2, j = a_{i} $ and
$ j' = 1 $.

\noindent
There is an arrow from the black vertex
$ V' $ to a white vertex
$ D_{i,j} $ if
$ i = h $ and
$ j = a_{h} $.
There are no other arrows.

We say that
$ D_{i,j} $ is an {\it even} (resp. {\it odd}) vertex if so is
$ i,  $ and
$ D_{i,j} $ is the {\it first} (resp. {\it last}) white vertex of
$ E(k,n) $ if
$ (i,j) = (1,1)  $ (resp. 
$ (h,a_{h})).\quad  $

Let
$ (k_{i}, n_{i}) $ be pairs of relatively prime positive integers such that
$ n_{i} > 1 \,(1 \le  i \le  g) $.
We define an oriented graph
$ E(k_{1},n_{1}; \dots ;k_{g},n_{g}) $ by induction on
$ g $.
It has a unique black vertex, and there is a unique white vertex, called the {\it last}
vertex, to which there is an arrow from the black vertex.
If
$ g = 1,  $ this is 
$ E(k_{1},n_{1}) $ defined above.
For
$ g > 1 $,
$ E(k_{1},n_{1}; \dots ;k_{g},n_{g}) $ is obtained by deleting the black vertex
of
$ E(k_{1},n_{1}; \dots ;k_{g-1},n_{g-1}) $ and the arrow between the black vertex
and the last white vertex, and then identifying the last white vertex of
$ E(k_{1},n_{1}; \dots ;k_{g-1},n_{g-1}) $ with the first vertex of
$ E(k_{g}+n_{g},n_{g}) $.

The white vertices of
$ E(k_{1},n_{1}; \dots ;k_{g},n_{g}) $ are naturally identified with the 
disjoint union of the white vertices of
$ E(k_{\nu },n_{\nu }) \,(1 \le  \nu  \le  g),  $ and the vertex 
corresponding to
$ D_{i,j} $ of
$ E(k_{\nu },n_{\nu }) $ will be denoted by
$ {D}_{i,j}^{(\nu )} $.

The following proposition is due to Deligne [4] and Zariski [20].

\proclaim{{\bf 2.6.~Proposition}}
Let
$ (V,0) $ be a germ of a reduced and irreducible plane curve, and
$ (k_{1}, n_{1}), \dots , (k_{g}, n_{g}) $ the Puiseux pairs of
$ (V,0) $.
Then the Enriques diagram of
$ (V,0) $ is
$ E(k_{1},n_{1}; \dots ;k_{g},n_{g}) $.
More precisely, the canonical embedded resolution
$ \phi  : (X', D) \rightarrow  ({\Bbb C}^{2},0) $ of
$ (V,0) $ is a composition of
$ \psi _{\nu } : (X_{\nu }, D_{\nu }) \rightarrow  (X_{\nu -1}, 
D_{\nu -1}) \,(1 \le  \nu  \le  g)  $ such that the proper transform
$ V_{\nu } $ of
$ V $ by
$ \phi _{\nu } = \psi _{1}\scirc \cdots \scirc \psi _{\nu } $ has 
modified Puiseux pairs
$ (k_{\nu +1}, n_{\nu +1}), \dots , (k_{g}, n_{g}) $ with respect to some local coordinates
$ x_{\nu }, y_{\nu } $ such that
$ x_{\nu } $ is the defining equation of the divisor corresponding to the
last white vertex of
$ E(k_{\nu },n_{\nu }),  $ where
$ (X_{g}, D_{g}) = (X', D), (X_{0}, D_{0}) = ({\Bbb C}^{2}, 0),  $ and
$ D_{\nu } = {\phi }_{\nu }^{-1}(0) $.
The resolution processes of
$ \psi _{\nu } $ and
$ \phi _{\nu } $ are expressed respectively by the oriented graphs
$ E(k_{\nu },n_{\nu }) $ and
$ E(k_{1},n_{1}; \dots ;k_{\nu },n_{\nu }) $.
\endproclaim

\noindent
{\it Remark.} For
$ \nu  < g,  $ the proper transform
$ V_{\nu } $ intersects transversally with the divisor corresponding
to the last white vertex of
$ E(k_{\nu },n_{\nu })  $ if
$ k_{\nu +1} > n_{\nu +1},  $ and they intersect maximally otherwise.
In both cases, their intersection number is equal to
$ n_{\nu +1} \cdots n_{g} $.

\medskip

This proposition is a direct consequence of the following inversion 
formula of Puiseux pairs due to Abhyanker [1], Deligne [4] and Zariski 
[19]:

\proclaim{{\bf 2.7.~Lemma}}
With the notation of (2.1), let
$ (k_{1}, n_{1}), (k_{2}, n_{2}), \dots , (k_{g}, n_{g}) $ be the modified Puiseux pairs
of
$ (V,0) $ with respect to coordinates
$ x, y $.
Then the modified Puiseux pairs of
$ (V,0) $ with respect to the coordinates
$ y, x $ are
$ (n_{1}, k_{1}), (k_{2}, n_{2}), \dots , (k_{g}, n_{g}) $.
\endproclaim

\noindent
{\bf 2.8.}
Let
$ (V,0) $ be as in (2.2).
Let
$ \Gamma  $ be the dual graph of the exceptional divisor
$ D $ of the canonical resolution
$ \phi  $ in (2.2).
By definition, the dual graph consists of white and black vertices and edges
connecting them such that:

\noindent
(i) The vertices are the same as the Enriques diagram.

\noindent
(ii) Two vertices are connected if the corresponding divisors intersect.

For relatively prime positive integers
$ (k, n) $ such that
$ n > 1,  $ we define a graph
$ \Gamma (k,n)  $ as follows:

\noindent
(i) The vertices are the same as
$ E(k,n) $.

\noindent
(ii) Two vertices
$ D_{i,j}, D_{i',j'} \,(i < i' $ or
$ i = i', j < j')  $ are connected in the following cases:

(a)
$ i' = i, j' = j + 1 $.

(b)
$ i' = i + 2, j = a_{i}, j' = 1 $.

(c)
$ i = h - 1, i' = h, j = a_{h-1}, j' = a_{h} $.

\noindent
The last white vertex and the black vertex are connected.
The vertices are not connected in the other cases.
(This means that, deleting the last white vertex, the white vertices 
have two connected components: one consists of even vertices and the 
other of odd ones, and they are connected linearly.)

We define a graph
$ \Gamma (k_{1},n_{1}; \dots ;k_{g},n_{g}) $ inductively by identifying 
the black vertex of
$ \Gamma (k_{1},n_{1}; \dots ;k_{g-1},n_{g-1}) $ 
with the last white vertex of 
$ \Gamma (1,n_{g}) $ if 
$ k_{g} = 1 $, and with the first odd vertex of
$ \Gamma (k_{g},n_{g}) $ otherwise, where the first odd vertex means 
the first vertex
$ D_{1,1} $ if
$ k_{g} > n_{g} $ and
$ D_{3,1} $ otherwise.
 (Note that the white vertices of 
$ \Gamma (1,n_{g}) $ are all even.)

Let
$ D_{\alpha }, D_{\beta } $ be irreducible components of
$ D $ in (2.2).
They intersect if and only if one of them is proximate to the other and there
are no vertices proximate to both.
So we get

\proclaim{{\bf 2.9.~Lemma}}
With the above notation, assume
$ (V,0) $ irreducible with Puiseux pairs
$ (k_{1}, n_{1}), (k_{2}, n_{2}), \dots , (k_{g}, n_{g}) $.
Then the dual graph
$ \Gamma  $ is
$ \Gamma (k_{1},n_{1}; \dots ;k_{g},n_{g}) $.
\endproclaim

\bigskip\bigskip
\centerline{\bf 3. Multiplicity}

\bigskip
\noindent
{\bf 3.1.}
With the notation of (2.2), let
$ m_{\alpha } $ be the multiplicity of
$ \phi ^{*}f $ along
$ D_{\alpha } (\alpha  \in  \Lambda ),  $ where
$ f $ is a reduced defining equation of
$ V $.
Let
$ P_{\alpha } $ be the center of the blow-up such that
$ D_{\alpha } $ is the proper transform of its exceptional divisor.
Let
$ V_{\alpha } $ be the proper transform of
$ V $ at
$ P_{\alpha } $.
Then we have
$$
m_{\alpha } = \sum _{\beta \leftarrow \alpha } m_{\beta } + 
 \text{mult} _{P_{\alpha }}V_{\alpha },
\leqno(3.1.1)
$$
where
$ \beta  \leftarrow  \alpha  $ if
$ D_{\alpha } $ is proximate to
$ D_{\beta } $.

\medskip
\noindent
{\bf 3.2.}
We study the 
$ \nu ^{th}  $ step of the resolution
$ \psi _{\nu } : X_{\nu } \rightarrow  X_{\nu -1} $ in the notation of (2.
6).
Let
$ k = k_{\nu }, n = n_{\nu } $ and
$ n' = n_{\nu +1} \cdots n_{g} $.
We have the continued fraction expansion (2.4.1) of
$ k/n $.
Let
$ r_{i} $ be as in (2.4.2).
Then we see

\medskip
\cond[3.2.1]
$ r_{i}n' $ is the multiplicity of the proper transform of 
$ V $ at the center of the blow up corresponding to 
$ {D}_{i,j}^{(\nu )} $.

\medskip\noindent
Let
$ {m}_{i,j}^{(\nu )}  $ be the multiplicity of
$ \phi ^{*}f $ along
$ {D}_{i,j}^{(\nu )} $ in the notation of (2.5).
Let
$ m_{i} = {m}_{i,{a}_{i}}^{(\nu )} $ and
$ \overline{m}_{i} = m_{i}/n' $ for
$ 1 \le  i \le  h $.
By (3.1.1) and (3.2.1), we get
$$
\overline{m}_{i} = a_{i}(r_{i} + \overline{m}_{i-1}) + \overline{m}_{i-2}\quad 
 \text{for\quad 1  } \le  i \le  h,
\leqno(3.2.2)
$$
where
$ \overline{m}_{0} = 0 $ and
$ \overline{m}_{-1} = m_{-1}/n' $ with
$ m_{-1} $ the multiplicity of
$ \phi ^{*}f $ along the divisor corresponding to the last white vertex
of
$ E(k_{\nu -1},n_{\nu -1}) $ if
$ \nu  > 1,  $ and
$ 0 $ otherwise.
We have
$ \overline{m}_{i} \in  {\Bbb Z} $ applying (3.2.1) and (3.2.2) 
inductively to
$ \nu ' < \nu  $.

\medskip
\noindent
{\it Remark.} If
$ k < n,  $ we have
$ a_{1} = 0, \overline{m}_{1} = \overline{m}_{-1} $ and
$ r_{2} = k $.
If
$ k = 1,  $ we have
$ h = 2, a_{1} = 0, a_{2} = n, r_{0} = r_{2} = 1 $ and
$ r_{1} = n $.

\proclaim{{\bf 3.3.~Lemma}}
Let 
$ k, n $ be as above, and
$ P_{1}, \dots , P_{h} $ and
$ Q_{1}, \dots , Q_{h} $ as in (2.4).
Then we have
$$
\overline{m}_{i}=
\cases
n{P}_{i}+{\overline{m}}_{-1}{Q}_{i}
&\text{if \,$ i : $ odd}
\\
({\overline{m}}_{-1}+k){Q}_{i}
&\text{if \,$ i : $ even.}
\endcases
\leqno(3.3.1)
$$
\endproclaim

\demo\nofrills {Proof.\usualspace}
We have
$$
r_{i+1} = - a_{i}r_{i} + r_{i-1}\quad  \text{for }1 \le  i \le  h
\leqno(3.3.2)
$$
by (2.4.2).
Together with (3.2.2), this implies
$$
(\overline{m}_{i} + r_{i+1}/2) = a_{i}(\overline{m}_{i-1} + r_{i}/2) + (\overline{m}
_{i-2} + r_{i-1}/2)\quad  \text{for\quad 1  } \le  i \le  h,
\leqno(3.3.3)
$$
where
$ \overline{m}_{-1} + r_{0}/2 = \overline{m}_{-1} + k/2, \overline{m}_{0} + r_{1}/2 
=n/2 $.
So we get
$$
\overline{m}_{i} + r_{i+1}/2 = (n/2)P_{i} + (\overline{m}_{-1} + k/2)Q_{i}\quad  \text{for }
1 \le  i \le  h
\leqno(3.3.4)
$$
On the other hand, we have
$$
(-1)^{i+1}r_{i+1} = - nP_{i} + kQ_{i}\,\,\, \text{for }1 \le  i \le  h
\leqno(3.3.5)
$$
by (3.3.2).
So we get the assertion.
\enddemo

\noindent
{\bf 3.4.}
With the notation of (3.2), let
$ C_{0} $ denote the last white vertex
$ {D}_{h,{a}_{h}}^{(\nu )} $ of
$ E(k_{\nu },n_{\nu }) $,
and
$ C_{1}, C_{2}, C_{3} $ the divisors intersecting with
$ C_{0} $ such that

\noindent
(i)
$ C_{1} $ is the last white vertex of
$ E(k_{\nu -1},n_{\nu -1}) $ if
$ k_{\nu } = 1,  $ and
$ {D}_{h-1,{a}_{h-1}}^{(\nu )} $ otherwise, 

\noindent
(ii)
$ C_{2} = {D}_{h,{a}_{h}-1}^{(\nu )}  $, 

\noindent
(iii)
$ C_{3} $ is the first odd white vertex of
$ E(k_{\nu +1},n_{\nu +1}) $ if
$ \nu  < g, k_{\nu +1} \ne  1  $, the last white vertex of
$ E(k_{\nu +1},n_{\nu +1}) $ if
$ \nu  < g, k_{\nu +1} = 1 $, and
$ V' $ otherwise.

Let
$ m'_{i} $ denote the multiplicity of
$ \phi ^{*}f $ along
$ C_{i},  $ and
$ \overline{m}'_{i} = m'_{i}/n',  $ where
$ n' = n_{\nu +1} \cdots n_{g}, k = k_{\nu }, n = n_{\nu },  $ and
$ \overline{m}_{-1} $ are as in (3.2).
Since the intersection number of
$ C_{0} $ with the total transform
$ \sum _{\alpha \in \Lambda } m_{\alpha }D_{\alpha } + V' $ of
$ V $ is zero, we have
$$
\sum _{i} \overline{m}'_{i} \equiv  0 \mod \overline{m}'_{0}.
\leqno(3.4.1)
$$
Since
$$
\overline{m}'_{2} = \overline{m}'_{0} - \overline{m}'_{1} - 1
\leqno(3.4.2)
$$
by (3.1.1) and (3.2.1), we get
$$
\overline{m}'_{3} \equiv  1 \mod \overline{m}'_{0}.
\leqno(3.4.3)
$$
Let
$ P = P_{h-1}, Q = Q_{h-1} $ in the notation of (3.3).
Then
$ kQ - nP = (-1)^{h} $ by (2.4).
So we get
$$
\align
\overline{m}'_{0} 
&= (\overline{m}_{-1} + k)n
\tag 3.4.4
\\
\overline{m}'_{1} 
&= \cases
({\overline{m}}_{-1}+k)Q-1 
&\text{if $ h $ : even}
\\
({\overline{m}}_{-1}+k)Q 
&\text{if $ h $ : odd.}
\endcases
\tag 3.4.5
\endalign
$$
by (3.3.1).
They imply
$$
\overline{m}'_{2} = \cases
({\overline{m}}_{-1}+k)(n-Q)
&\text{if $ h $ : even}
\\
({\overline{m}}_{-1}+k)(n-Q)-1 
&\text{if $ h $ : odd.}
\endcases
\leqno(3.4.6)
$$
by (3.4.2).
We can verify that
$ \overline{m}'_{3} $ is equal to
$ 1 $ if
$ \nu  = g, 1 + \overline{m}'_{0} $ if
$ \nu  < g, k_{\nu +1} = 1,  $ and
$ 1 + \overline{m}'_{0}(1 + [n_{\nu +1}/k_{\nu +1}]) $ otherwise.
But it will not be used later, because (3.4.3) is sufficient.

Let
$ w_{\nu } $ be as in (1.4).
Then (3.4.4) implies
$$
w_{\nu } = \overline{m}_{-1} + k_{\nu },\,\,\,m'_{0} = w_{\nu }n_{\nu } 
\cdots  n_{g}
\leqno(3.4.7)
$$
by induction on
$ \nu  $.
We can also verify that the multiplicity on the first even vertex of
$ E(k_{\nu },n_{\nu }) $ is
$ w_{\nu }n_{\nu +1} \cdots n_{g} $.

\medskip\noindent
{\it Remark.} It is known by a topological method that the multiplicity
$ {m}_{h,{a}_{h}}^{(\nu )}  $ of
$ \phi ^{*}f $ along the last white vertex
$ {D}_{h,{a}_{h}}^{(\nu )} $ of
$ E(k_{\nu },n_{\nu }) $ is
$ w_{\nu }n_{\nu } \cdots  n_{g} $.
See e.g. [6].
(Note that the last white vertices correspond to the rapture points of 
the resolution graph, and hence to the splice components of the 
$ EN $-diagram, and the multiplicity can be interpreted as the linking 
number for each splice component.)
The assertions (3.4.5--6) imply that if we choose integers
$ \beta $,
$ \beta ' $ such that
$ \beta w_{\nu } \equiv  1 \mod  n_{\nu } $ and
$ \beta 'n_{\nu } \equiv  1  \mod  w_{\nu } $,
then the multiplicities along the adjacent divisors
$ C_{1} $,
$ C_{2} $ are given modulo
$ {m}_{h,{a}_{h}}^{(\nu )} $ by
$ - \beta w_{\nu }n_{\nu +1} \cdots  n_{g} $ and
$ - \beta 'n_{\nu } \cdots  n_{g} $ respectively if
$ h $ is odd, and the order is reversed if
$ h $ is even.
This assertion in a more general situation is remarked in [12], p. 127 
without a reference.
It is clear that this property about the multiplicities modulo
$ {m}_{h,{a}_{h}}^{(\nu )} $ is enough to show (1.5).
See also [9].
(Note that the argument in 3.1 of loc. cit. is slightly misstated because 
any numbers
$ s $ in
$ [1,m_{w_{i}}-1] $ cannot always be written as stated there.)

\medskip
\noindent
{\bf 3.5.}~{\it Remark.} With the notation of (2.2), let
$ d_{\alpha }  $ be the multiplicity of the determinant of
$ d\phi  $ along
$ D_{\alpha },  $ and
$ \tilde{d}_{\alpha } = d_{\alpha } + 1, e_{\alpha } = \tilde{d}_{\alpha }/m
_{\alpha } $.
Assume
$ V $ irreducible.
Let
$ {d}_{i,j}^{(\nu )} = d_{\alpha }, {\tilde{d}}_{i,j}^{(\nu )} = d_{\alpha } $ and
$ {e}_{i,j}^{(\nu )} = e_{\alpha },  $ if
$ D_{\alpha } = {D}_{i,j}^{(\nu )} $.
We say that
$ {D}_{i,j}^{(\nu )} $ is odd (resp. even) if
$ i $ is odd (resp. even) or
$ {D}_{i,j}^{(\nu )} $ is the last white vertex of
$ E(k_{\nu },n_{\nu }) $.
Then, for
$ {D}_{i,j}^{(\nu )} \ne  {D}_{i',j'}^{(\nu ')},  $ we have the inequality
$ {e}_{i,j}^{(\nu )} > {e}_{i',j'}^{(\nu ')} $ in the following cases:

(a)
$ {D}_{i,j}^{(\nu )}, {D}_{i',j'}^{(\nu ')} $ are odd, and
$ (\nu , i, j) > (\nu ', i', j') $ with
$ \nu ' > 1 $ or
$ (i, j) < (i', j') $ with
$ \nu  = \nu ' = 1 $.

(b)
$ {D}_{i,j}^{(\nu )}, {D}_{i',j'}^{(\nu ')} $ are even,
$ \nu  = \nu ' $ and
$ (i, j) < (i', j') $.

(c)
$ {D}_{i,j}^{(\nu )} $ is even and
$ {D}_{i',j'}^{(\nu ')} $ is odd, 
$ \nu  \ge  \nu ' > 1 $.

(d)
$ {D}_{i',j'}^{(\nu ')} $ is the last white vertex of
$ E(k_{1},n_{1}) $.

\noindent
Here
$ > $ denote the lexicographic order.
The case 
$ (c)  $ follows from 
$ (a)  $ and 
$ (b),  $ and may be omitted.
By 
$ (d),  $ the minimum of
$ {e}_{i,j}^{(\nu )} $ is attained by the last white vertex of
$ E(k_{1},n_{1}) $.
We can show that it is equal to
$ (k_{1}+n_{1})/k_{1}n_{1} \cdots n_{g} $.

In fact, with the notation of (3.2), let
$ \tilde{d}_{i} = d_{i} + 1 \,(-1 \le  i \le  h) $ with
$ d_{i} = {\tilde{d}}_{i,{a}_{i}}^{(\nu )} \,(1 \le  i \le  h) $,
$ d_{0} = 0 $ and
$ d_{-1} $ the multiplicity of
$ \det d\phi  $ along the divisor corresponding to the last white vertex
of
$ E(k_{\nu -1},n_{\nu -1}) $ if
$ \nu  > 1,  $ and
$ 0 $ otherwise.
Then we have
$$
{\tilde{d}}_{i,j}^{(\nu )} = j\tilde{d}_{i-1} + \tilde{d}_{i-2},\,\,\,\tilde{d}_{i} = a_{i}\tilde{d}
_{i-1} + \tilde{d}_{i-2}\quad  \text{for }1 \le  i \le  h,
$$
and
$ \tilde{d}_{i} = P_{i} + \tilde{d}_{-1}Q_{i} $ by (2.4.3).
Combining with (2.4.4), we can verify the assertion, 
showing also
$ {e}_{i,j}^{(\nu )} < 1 $ inductively, which implies
$ n\tilde{d}_{-1}/\overline{m}_{-1} < 1 $ in the above notation.

\bigskip\bigskip
\centerline{\bf 4. Proof of Theorem (1.5)}

\bigskip
\noindent
{\bf 4.1.}
With the notation of (1.1), let
$$
{H}_{\lambda }^{p,q} = {\Gr}_{F}^{p}{\Gr}_{p+q}^{W}H^{n}
(X_{\infty },{\Bbb C})_{\lambda }.
\leqno(4.1.1)
$$
where
$ F $ and
$ W $ are the Hodge and weight filtrations of the mixed Hodge 
structure.
We denote by
$ D_{\beta } \,(\beta  \in  \Lambda ') $ the irreducible components of 
the proper transform
$ V' $ of
$ V (cf. (2.2)) $.
Let
$ \Lambda _{\alpha } = \{\beta  \in  \Lambda  \cup  \Lambda ' : D
_{\alpha } \cap  D_{\beta } \ne  \emptyset \},  $ and
$ m_{\beta } = 1 $ for
$ \beta  \in  \Lambda ' $.
By [15, (3.13--14)] (see also [13]), we have
$$
\aligned
{H}_{\lambda }^{0,1}
&= \bigoplus 
\Sb \alpha \in \Lambda\\
0\le c\le m_{\alpha }\\
\lambda =\bold{e}(c/m_{\alpha})
\endSb
H(\alpha ,c)\quad \text{with}
\\
H(\alpha ,c)
&= H^{1}\biggl (\bold{P}^{1}, {\Cal O}_{\bold{P}^{1}}
\biggl(\sum _{\beta \in \Lambda _{\alpha }}
\biggl(-{cm_{\beta }\over m_{\alpha }} +
\biggl[{cm_{\beta }\over m_{\alpha }}\biggr]\biggr)\biggr)\biggr).
\\
&\simeq
H^{0}\biggl(\bold{P}^{1}, {\Cal O}_{\bold{P}^{1}}
\biggl(\sum _{\beta \in \Lambda _{\alpha }}
\biggl({cm_{\beta }\over m_{\alpha }} - 
\biggl[{cm_{\beta }\over m_{\alpha }}\biggr]\biggr) - 2\biggr)\biggr).
\endaligned
\leqno(4.1.2)
$$
We have
$ H(\alpha ,c) = 0 $ if
$ |\Lambda _{\alpha }| \le  2,  $ and
$ H(\alpha ,c) = {\Bbb C} $ if and only if
$$
\sum _{\beta \in \Lambda _{\alpha }}
\left({cm_{\beta }\over m_{\alpha }} - 
\left[{cm_{\beta }\over m_{\alpha }}\right]\right)  = 2.
\leqno(4.1.3)
$$
{\it Remark.} We can show (4.1.2) using the weight spectral sequence of
$ H^{\ssbull }(X_{\infty },{\Bbb C}) $.
In fact, let
$ g = \phi ^{*}f $.
Then we can calculate
$ {\Gr}_{0}^{W}\psi _{g,\lambda }{\Bbb C}_{X'},  $ where
$ \psi _{g,\lambda }{\Bbb C}_{X'} $ is the 
$ \lambda  $-eigenvalue part of the nearby cycle
$ \psi _{g,\lambda }{\Bbb C}_{X'} $ by the action of the monodromy, and
$ W $ is the monodromy filtration.
Its direct factor supported on
$ D_{\alpha } $ is an intersection complex associated with a local 
system of rank one whose monodromy around the intersection with
$ D_{\beta } $ is
$ \lambda ^{-m_{\beta }} $.
Then we can deduce (4.1.2) using Deligne's canonical extension of the local
system.

\medskip

For the calculation of (4.1.2) we use the following

\proclaim{{\bf 4.2.~Lemma}}
Let 
$ \alpha _{1}, \alpha _{2}, \alpha _{3} $ be three real numbers such that
$ \alpha _{1} + \alpha _{2} + \alpha _{3} \in  {\Bbb Z} $.
Then
$ \sum _{i} (\alpha _{i} - [\alpha _{i}]) = 2 $ if and only if
$ \alpha _{1} \notin  {\Bbb Z} $ and
$ \alpha _{2} + \alpha _{3}  - [\alpha  _{2} + \alpha _{3} ] 
< \alpha _{3} - [\alpha _{3}] $.
\endproclaim

\noindent
{\bf 4.3.}~{\it Proof of} (1.5).
Since
$ V $ is irreducible, the monodromy is semisimple and one is not an eigenvalue
of the monodromy [2], [8], so that
$$
{H}_{\lambda }^{0,0} = {H}_{\lambda }^{1,1} = {H}_{1}^{p,q} = 0,
\leqno(4.3.1)
$$
cf. [15].
So it is enough to calculate
$ {H}_{\lambda }^{0,1} $ for
$ \lambda  \ne  1 $ by the Hodge symmetry.
We have
$ H(\alpha ,c) = 0 $ unless
$ D_{\alpha } $ is the last white vertex of
$ E(k_{\nu },n_{\nu }) $ for some
$ \nu , cf. (2.9) $.
So we may assume that
$ D_{\alpha } $ is
$ C_{0} $ in (3.4).
The multiplicities are calculated in (3.4).
Let
$ w = w_{\nu } $.
By (3.4.7) applied to
$ \nu  - 1,  $ we get
$ n | \overline{m}_{-1},  $ and
$ (w, n) = 1 $.
With the notation of (3.2), the integer
$ c $ in
$ [0, wnn') $ is expressed uniquely as
$ c = iw + jn + rwn $ for
$ i, j, r \in  {\Bbb Z} $ such that
$$
0 \le  i < n, 0 \le  iw + jn < wn, 0 \le  r < n'.
$$
By (3.4.5--6), (4.1.3) and (4.2), 
$ H(\alpha ,c) = {\Bbb C} $ if and only if
$$
\alignat2
&\frac{c(n - Q)}n \notin  {\Bbb Z}, \quad
\frac{cQ}n - \biggl[\frac{cQ}n\biggr] < 
\frac{c}{wn} - \biggl[\frac{c}{wn}\biggr]
&&\qquad\text{if  $ h $ : even}
\\
&\frac{cQ}n \notin  {\Bbb Z}, \quad
\frac{c(n - Q)}n - \biggl[\frac{c(n - Q)}n\biggr] < 
\frac{c}{wn} - \biggl[\frac{c}{wn}\biggr]
&&\qquad\text{if  $ h $ : odd.}
\endalignat
$$
Since
$ Qw \equiv  (-1)^{h}  \mod  n $ by (2.4.4--5), 
this condition is equivalent to
$$
i > 0, \,\, j > 0.
$$
So the contribution of
$ C_{0} $ to the exponents which are greater than 
$ 1 $ is given by
$$
\left\{{2-\left({{1 \over {n}_{\nu +1}\cdots
{n}_{g}}\left({{i \over {n}_{\nu }}+{j \over {w}_{\nu }}}\right)+
{r \over {n}_{\nu +1}\cdots {n}_{g}}}\right)}\right\}
$$
for
$ 0 < i < n_{\nu } $,
$ 0 < j < w_{\nu } $,
$ 0 \le  r < n_{\nu +1} \cdots n_{g} $ such that
$ i/n_{\nu } + j/w_{\nu } < 1 $.
See (1.2) and (4.1).
Then the assertion follows from the symmetry of exponents, 
cf. (1.2.3).

\bigskip\bigskip
\centerline{{\bf 5. Variance of exponents}}

\bigskip
\noindent
In this section we prove the following conjecture of C. Hertling [7] 
in the case of irreducible plane curve singularities.

\bigskip
\noindent
{\bf 5.1.~Conjecture.} 
$ V_{f} := {1\over \mu}\sum _{i} (\alpha _{i} - {n+1 \over 2})^{2} 
\le (\alpha _{\mu } - \alpha _{1})/12 $ with the notation of (1.2). 

\medskip
\noindent
{\it Remark.} 
$ V_{f} $ is the variance (i.e. the square of the standard deviation) of 
the exponents.
Hertling (loc. cit.) showed that the equality holds in (5.1) if
$ f $ is quasihomogeneous.
In particular, (5.1) is true if 
$ g = 1 $ in the case of irreducible plane curve singularity.

\proclaim{{\bf 5.2.~Theorem}}
$ V_{f} < (\alpha _{\mu } - \alpha _{1})/12 $ if 
$ g > 1 $ in the case of irreducible plane curve singularity.
\endproclaim

\demo\nofrills {Proof.\usualspace}
Let
$ \{\alpha _{i}\}_{i\in \Lambda (\nu )} $ be the exponents which are 
less than
$ 1 $ and come from the
$ \nu $th  part of the Enriques diagram corresponding to
$ E(k_{\nu },n_{\nu }) $, cf. (2.5).
Here
$ \{\Lambda (\nu )\}_{1 \le  \nu  \le  g} $ is a partition of
$ \{1, \dots , \mu /2\} $.
By (1.5),
$ \{\alpha _{i}\}_{i\in \Lambda (\nu )} $ are given by (1.5.2) with
$ \nu  $ fixed.

Let
$ n'_{\nu } = n_{\nu +1} \cdots  n_{g} $ for
$ 0 \le  \nu  \le  g $ (where
$ n'_{g} = 1) $.
We  define
$$
\aligned
\mu ^{(\nu )} 
&= 2|\Lambda (\nu )| \,(= (w_{\nu } - 1)(n_{\nu } - 1)n'_{\nu }),
\\
S^{(\nu )} 
&= 2\sum _{i\in \Lambda (\nu )} (\alpha _{i} - 1)^{2},
\\
\varepsilon ^{(\nu )} 
&= 6S^{(\nu )} - \mu ^{(\nu )}(1 - \alpha _{1}).
\endaligned
$$
Then we have to show
$ \sum_{\nu =1}^{g} \varepsilon ^{(\nu) } < 0 $.

By (1.5), we have
$$
\alpha _{1} = {k_{1} + n_{1}\over k_{1}n'_{0}},
\leqno(5.2.1)
$$
because
$ (k_{1} + n_{1})/k_{1}n_{1} < 1 $.
We can verify (see also [7])
$$
6S^{(g)} = \mu ^{(g)}\left(1 - {1\over w_{g}} - {1\over n_{g}}\right).
\leqno(5.2.2)
$$

Let
$ a = w_{\nu } $,
$ b = n_{\nu } $,
$ c = n'_{\nu } $,
and
$ \Lambda (a,b) = \left\{(i,j) \in  {\Bbb Z}^{2} : i, j > 0, 
{i\over a} + {j\over b} \le 1\right\} $.
Then
$$
\aligned
6S^{(\nu )} = 
&\,\sum_{k=0}^{c-1} \,\sum_{(i,j)\in \Lambda (a,b)}
{12\over c^{2}}\left(1 - {i\over a} - {j\over b} + k\right)^{2}
\\
=
&\,\sum_{(i,j)\in \Lambda (a,b)} \left({12\over c}
\left(1 - {i\over a} - {j\over b}\right)^{2} + 
12\left(1 - {i\over a}- {j\over b}\right){c-1\over c}\right)
\\
& \quad + (a - 1)(b - 1){(c - 1)(2c - 1)\over c}
\\
\le  
&\,{(a - 1)(b - 1)\over c}\left(\left(1 - {1\over a} - {1\over b}\right) 
+ 2(c - 1) + (c - 1)(2c - 1)\right)
\\
& \quad - \left(1 - {1\over b}\right)\left(1 - {1\over c}\right)
(a + b - 1)
\endaligned
$$
by (5.3) below together with a formula similar to (5.2.2) (where
$ w_{g} = a $,
$ n_{g} = b) $.
So we get
$$
\varepsilon ^{(\nu )} \le  (a - 1)(b - 1)c - ab + {b - 1\over ac} 
+ {a - 1\over b} + 1 + (a - 1)(b- 1)c\alpha _{1}.
$$
Since
$ w_{\nu } = w_{\nu -1}n_{\nu -1}n_{\nu } + k_{\nu } $,
this implies for
$ 1 < \nu  < g $
$$
\aligned
\varepsilon ^{(\nu )} 
&\le  w_{\nu }n'_{\nu -1} - w_{\nu -1}n'_{\nu -2} - w_{\nu }n_{\nu } 
+ w_{\nu -1}n_{\nu -1} - (k_{\nu } - 1)
\left( n'_{\nu } - {1\over n_{\nu }}\right)
\\
& \quad - (n'_{\nu -1} - 1) + {n_{\nu } - 1\over w_{\nu }n'_{\nu }}+ 
(w_{\nu } - 1)(n_{\nu } - 1)n'_{\nu }\alpha _{1}.
\endaligned
$$
Then, after a calculation, we get a simple expression
$$
\sum_{\nu =1}^{g} \varepsilon ^{(\nu )} \le  \sum_{\nu =1}^{g}
\left({n_{\nu } - 1\over w_{\nu }n'_{\nu }} - 
(k_{\nu } - 1)\left(n'_{\nu } - {1\over n_{\nu }}\right) 
- (n'_{\nu -1} - 1)\right) + \mu \alpha _{1}.
\leqno(5.2.3)
$$

By induction on
$ j $,
we see
$$
\sum_{\nu =1}^{j} \mu ^{(\nu )} = (w_{j} - 1)n'_{j-1} - \sum_{\nu =1}^{j}
(k_{\nu } + n_{\nu -1}n_{\nu } - 1)n'_{\nu },
\leqno(5.2.4)
$$
where
$ n_{0} = 0 $.
Combined with
$ w_{g}n_{g} = {\sum\nolimits}_{\nu =1}^{g} 
k_{\nu }n'_{\nu -1}n'_{\nu } $, this (for
$ j = g $) implies
$$
\aligned
\mu \alpha _{1} - \sum_{\nu =1}^{g} (k_{\nu } - 1)
\left(n'_{\nu } - {1\over n_{\nu }}\right) 
&= \sum_{\nu =1}^{g} (k_{\nu } - 1)
\left(n'_{\nu } - {1\over n_{\nu }}\right)(n'_{\nu -1}\alpha _{1} - 1)
\\
& \qquad +\sum_{\nu =1}^{g} (n'_{\nu } - 1)n'_{\nu -1}\alpha _{1}.
\endaligned
$$
Since
$ n'_{\nu -1}\alpha _{1} - 1 < 0 $ for
$ \nu  > 1 $ and
$ n'_{0}\alpha _{1} - 1 = n_{1}/k_{1} $,
it remains to show
$$
(k_{1} - 1)\left(n'_{1} - {1\over n_{1}}\right)
{n_{1}\over k_{1}} + (n'_{1} - 1){n_{1}\over k_{1}} + 
\sum_{\nu =1}^{g} {n_{\nu } - 1\over w_{\nu }n'_{\nu }}
< n'_{0} - 1.
$$
But this is reduced to
$ {\sum\nolimits}_{\nu =1}^{g} (n_{\nu } -1)/w_{\nu }n'_{\nu } 
< (n_{1} - 1)/k_{1} $,
and is easily verified.
\enddemo

\proclaim{{\bf 5.3.~Lemma}}
Let 
$ a, b $ be relatively prime positive integers such that
$ a > b $.
Then
$$
\sum_{(i,j)\in \Lambda (a,b)} \left(1 - {i\over a} - {j\over b}\right) 
\le  {(a - 1)(b - 1)\over 6} - {(b - 1)(a + b - 1)\over 12b }.
\leqno(5.3.1)
$$
\endproclaim

\demo\nofrills {Proof.\usualspace}
Let
$ F(a,b) = \sum_{(i,j)\in \Lambda (a,b)} 
\left(1 - {i\over a} - {j\over b}\right) $.
We see
$$
\aligned
F(a,b) - {a-b\over a}F(a-b,b) 
&= \sum_{(i,j)\in \Lambda (b,b)} 
\left(1 - {i\over a} - {j\over b}\right)
\\
&={(b - 1)(2ab - b^{2} - a - b)\over 6a}.
\endaligned
\leqno(5.3.2)
$$
Let
$ a = mb + k $ with
$ 0 < k < b $.
Then (5.3.2) implies
$$
F(a,b) - {k\over a}F(k,b) = {(b-1)m\over 12a}
((2b^{2}-b)m + 4bk -2k - 3b).
\leqno(5.3.3)
$$
Let
$ E(a,b) = (a-1)(b-1)/6 - F(a,b) $.
Using (5.3.3) we can verify
$$
E(a,b) - {k\over a}E(b,k) = {m(b-1)(a+b+k)\over 12a}.
\leqno(5.3.4)
$$
So the inequality
$ E(a,b) \ge  (b - 1)(a + b - 1)/12b $ follows by induction, because
$$
{(k-1)(b+k-1)\over a} + {m(b-1)(a+b+k)\over a} \ge  
{(b-1)(a+b-1)\over b}.
\leqno(5.3.5)
$$
This completes the proof of (5.2).
\enddemo

\bigskip\bigskip
\centerline{{\bf References}}

\bigskip

\item{[1]}
S. Abhyanker, Inversion and invariance of characteristic pairs, Am. J.
Math. 89 (1967) 363--372.

\item{[2]}
N. A'Campo, Sur la monodromie des singularit\'es isol\'ees
d'hypersurfaces complexes, Inv. Math. 20 (1973) 147--169.

\item{[3]}
P. Deligne, Intersections sur les surfaces r\'eguli\`eres,
in SGA7 X, Lect. Notes in Math. vol. 340, Springer, Berlin, 1973, 
pp. 1--37.

\item{[4]}
P. Deligne, manuscript.

\item{[5]}
P. Deligne, Th\'eorie de Hodge I, Actes Congr\`es Intern. Math., (1970) 
425--430; II, Publ. Math. IHES, 40 (1971) 5--58; III, ibid. 44 (1974) 
5--77.

\item{[6]}
D. Eisenbud and W. Neumann, Three-Dimensional Link Theory and 
Invariants of Plane Curve Singularities, Ann. Math. Studies 110, 
Princeton Univ. Press, 1985.

\item{[7]}
C. Hertling, Variance of the spectral numbers, preprint (math.CV/0007187).

\item{[8]}
D.T. L\^e, Sur les noeuds alg\'ebriques, 
Compositio. Math. 25 (1972) 281--321.

\item{[9]}
A. N\'emethi, On the spectrum of curve singularities, in Proceedings 
of the Singularity Conference, Oberwolfach, Progr. in Math., 162, 
Birkh\"auser 1998, pp. 93--102.

\item{[10]}
M. Saito, On the exponents and the geometric genus of an isolated 
hypersurface singularity, Proc. Symp. Pure Math. 40 (1983) Part 2, 
465--472.

\item{[11]}
M. Saito, Exponents and Newton polyhedra of isolated hypersurface 
singularities, Math. Ann. 281 (1988) 411--417.

\item{[12]}
R. Schrauwen, Topological series of isolated plane curve singularities, 
Enseign. Math. 36 (1991), 115--141.

\item{[13]}
R. Schrauwen, J. Steenbrink and J. Stevens, Spectral pairs and topology 
of curve singularities, Proc. Symp. Pure Math. 53 (1991), 305--328.

\item{[14]}
J. Steenbrink, Intersection form for quasi-homogeneous singularities,
Compo. Math. 34 (1977) 211--223.

\item{[15]}
J. Steenbrink, Mixed Hodge structure on the vanishing cohomology, in
Real and Complex Singularities (Proc. Nordic Summer school, Oslo, 
1976, pp. 525--563) Alphen a/d Rijn: Sijthoff \& Noordhoff 1977.

\item{[16]}
A. Varchenko, Asymptotic Hodge structure in the vanishing cohomology,
Math. USSR Izv. 18 (1982) 465--512.

\item{[17]}
A. Varchenko, The asymptotics of holomorphic forms determine a mixed
Hodge structure, Soviet Math. Dokl., 22 (1980) 772--775.

\item{[18]}
A. Varchenko, The complex exponent of a singularity does not change 
along strata 
$ \mu $  = const, Func. Anal. Appl. 16 (1982) 1--9.

\item{[19]}
O. Zariski, Studies in equisingularity III, saturation of local rings
and equisingularity, Am. J. Math. 90 (1968) 961--1023.

\item{[20]}
O. Zariski, General theory of saturation and of saturated local ring
II, Am. J. Math. 93 (1971) 872--964.

\bigskip\noindent
Sept. 26, 2000

\bye